\documentclass[a4paper,openright,11pt]{article}
\usepackage[utf8]{inputenc} %utfx8
\usepackage{amsthm,amssymb,amsbsy,amsmath,amsfonts,amssymb,amscd,mathrsfs}
\usepackage{graphicx} % [dvips] se vuoi gli eps!
\usepackage{subfigure}
\usepackage{fancyhdr}
\usepackage[all]{xy}

\DeclareGraphicsExtensions{.jpg, .png}
\usepackage{caption}

\theoremstyle{plain}
\newtheorem{theorem}{Theorem}%[section]
\newtheorem{corollary}[theorem]{Corollary}
\newtheorem{proposition}[theorem]{Proposition}
\theoremstyle{definition}
\newtheorem{definition}[theorem]{Definition}
\newtheorem{lemma}[theorem]{Lemma}
\theoremstyle{remark}
\newtheorem{remark}[theorem]{Remark}
\newcommand{\bt}{\begin{theorem}}
\newcommand{\et}{\end{theorem}}
\newcommand{\bl}{\begin{lemma}}
\newcommand{\el}{\end{lemma}}
\newcommand{\bp}{\begin{proposition}}
\newcommand{\ep}{\end{proposition}}
\newcommand{\bc}{\begin{corollary}}
\newcommand{\ec}{\end{corollary}}
\newcommand{\bdeff}{\begin{definition}}
\newcommand{\edeff}{\end{definition}}
\newcommand{\brem}{\begin{remark}}
\newcommand{\erem}{\end{remark}}
\newcommand{\bi}{\begin{itemize}}
\newcommand{\iii}{\item}
\newcommand{\ei}{\end{itemize}}
\newcommand{\bd}{\begin{description}}
\newcommand{\ed}{\end{description}}

\newcommand{\bqn}{\begin{eqnarray}}
\newcommand{\eqn}{\end{eqnarray}}
\newcommand{\eqnn}{\nonumber\end{eqnarray}}

\newcommand{\ba}[1]{\begin{array}{#1}}
\newcommand{\ea}{\end{array}}
\newcommand{\distr}{\Delta}
\newcommand{\metr}{\mathbf{g}}

\newcommand{\lam}{\lambda}
\newcommand{\g}{\gamma}

\newcommand{\ph}{\varphi}

\newcommand{\R}{\mathbb{R}}
\newcommand{\N}{\mathbb{N}}

\newcommand{\mb}[1]{\mathbb{ #1 }}

\newcommand{\all}{\forall}

\newcommand{\la}{\langle}
\newcommand{\ra}{\rangle}
\newcommand{\ffoot}[1]{\footnote{{\noindent#1}}}
\newcommand{\VecM}{\mathrm{Vec}(M)}
\newcommand{\virg}[1]{``#1''}
\newcommand{\tx}[1]{\mathrm{#1}}
\newcommand{\til}[1]{\widetilde{#1}}

\begin{document}
\begin{center} \noindent
{\LARGE{\sl{\bf Sub-Riemannian structures on 3D Lie groups}}}
\vskip 0.6 cm
Andrei Agrachev\\ 
{\footnotesize SISSA, Trieste, Italy and MIAN, Moscow, Russia - {\tt agrachev@sissa.it}}\\
\vskip 0.3 cm
Davide Barilari\\ 
{\footnotesize SISSA, Trieste, Italy - {\tt barilari@sissa.it}}\\
\vskip 0.6cm
%%Ugo Boscain\symbolfootnote[0]{This research has been supported by the European Research Council, ERC StG 2009 GeCoMethods, contract number 239748, by the ANR GCM, and by the DIGITEO project CONGEO.}\\
%%{\footnotesize CNRS, CMAP Ecole Polytechnique, Paris, France - {\tt boscain@cmap.polytechnique.fr}}
%%\vskip 0.3cm
\today
\end{center}

\begin{abstract}
We give the complete classification of left-invariant sub-Riemannian
structures on three dimensional Lie groups in terms of the basic
differential invariants introduced in \cite{AAAICM,agrexp}. This classifications recovers other known classification results in the literature, in particular the one obtained in \cite{falbel} in terms of curvature invariants of a canonical connection. Moreover, we explicitly find a
sub-Riemannian isometry between the nonisomorphic Lie groups $SL(2)$
and $A^+(\R)\times S^1$, where $A^+(\R)$ denotes the group of
orientation preserving affine maps on the real line.\\
%This isometry could be also interpreted as special kind of coordinates on $SL(2)$.
\end{abstract}
%\tableofcontents

{\bf MSC classes}:	53C17 $\cdot$	22E30 $\cdot$ 49J15\\

{\bf Keywords}: sub-Riemannian geometry, Lie groups, left-invariant structures.
\section{Introduction}
In this paper, by a sub-Riemannian manifold we mean a triple $(M,\distr,\metr)$, where $M$ is a connected smooth manifold of dimension $n$, $\distr$ is a smooth vector distribution of constant rank $k<n$, and $\metr$ is a Riemannian metric on $\distr$, smoothly depending on the point.

In the following we always assume that the distribution $\distr$ satisfies the \emph{bracket generating condition} (also known as \emph{H$\ddot{o}$rmander condition}), i.e. the Lie algebra generated by vector fields tangent to the distribution spans at every point the tangent space to the manifold.

Under this assumption, $M$ is endowed with a natural structure of metric space, where the distance is the so called \emph{Carnot-Caratheodory distance}
\begin{eqnarray}
d(p,q)=
\inf\{\int_0^T\sqrt{\metr_{\g(t)}(\dot\gamma(t),\dot\gamma(t))}\,
dt~|~ \gamma:[0,T]\to M \mbox{ is a Lipschitz curve},\cr
\gamma(0)=p,\gamma(T)=q, ~~\dot \gamma(t)\in\distr_{\gamma(t)}\mbox{ a.e. in $[0,T]$} \}.\nonumber
\end{eqnarray}
As a consequence of the H\"ormander condition this distance is always finite and continuous, and induces on $M$ the original topology (see Chow-Rashevsky Theorem, \cite{agrcontrol}). Standard references on sub-Riemannian geometry are \cite{bellaiche,gromov,montgomery}.
%Under this assumption it follows from the classical Chowtheorem that $d$ is a well defined metric on $M$ and it induces on $M$ the original topology.

A sub-Riemannian structure is said to be \emph{contact} if its
distribution is defined as the kernel of a contact differential
one form $\omega$, i.e. $n=2m+1$ and
$\left(\bigwedge^md\omega\right)\wedge\omega$ is a nonvanishing
$n$-form on $M$.

In this paper we focus on the three dimensional case. Three
dimensional contact sub-Riemannian structures have been deeply studied in the last years (for example see \cite{miovolume,agrexp,localagr}) and they have two basic
differential invariants $\chi$ and $\kappa$ (see Section \ref{sec:inv} for the
precise definition and \cite{nostrolibro, agrexp} for their role in the
asymptotic expansion of the sub-Riemannian exponential map). 

The invariants $\chi$ and $\kappa$ are smooth real functions on
$M$. It is easy to understand, at least heuristically, why it is
natural to expect exactly {\it two} functional invariants. Indeed,
in local coordinates the sub-Riemannian structure is defined
by its orthonormal frame, i.e. by a couple of smooth vector
fields on $\mathbb R^3$ or, in other words, by 6 scalar functions
on $\mathbb R^3$. One function can be normalized by the rotation
of the frame within its linear hull and three more functions by
 smooth change of variables. What remains are two scalar
functions.

In this paper we exploit these local invariants to provide a
complete classification of left-invariant structures on 3D Lie
groups. A sub-Riemannian structure on a Lie group is said to be
\emph{left-invariant} if its distribution and the inner product are
preserved by left translations on the group. A left-invariant
distribution is uniquely determined by a two dimensional subspace
of the Lie algebra of the group. The distribution is bracket
generating (and contact) if and only if the subspace is not a Lie
subalgebra. %In particular, since we always assume the structure to be bracket generating,  

Left-invariant structures on Lie groups are the basic models of sub-Riemannian manifolds and the study of such structures is the starting point to understand the general properties of sub-Riemannian geometry. 
In particular, thanks to the group structure, in some of these cases it is also possible to compute explicitly the sub-Riemannian distance and geodesics (see in particular \cite{gersh} for the Heisenberg group, \cite{boscainrossi} for semisimple Lie groups with Killing form and \cite{sachkovmois,sachkov} for a detailed study of the sub-Riemannian structure on the group of motions of a plane). 

\emph{Remark.\ }
The problem of equivalence for several geometric structures close to left-invariant sub-Riemannian structures on 3D Lie groups were studied in several publications (see \cite{cartan2, cartan1, falbel, stric,stric2}). In particular in \cite{stric} the author provide a first classification of symmetric sub-Riemannian structures of dimension 3, while in \cite{falbel} is presented a complete classification of sub-Riemannian homogeneous spaces (i.e., sub-Riemannian structures which admits a transitive Lie group of isometries acting smoothly on the manifold) by means of an adapted connection. 
The principal invariants used there, denoted by $\tau_{0}$ and $K$, coincide, up to a normalization factor, with our differential invariants $\chi$ and $\kappa$.
\medskip

A standard result on the classification of 3D Lie algebras (see, for instance, \cite{jacob}) reduce the analysis on the Lie
algebras of  the following Lie groups:
\bi
\iii[] $H_{3}$,  the
Heisenberg group, 
\iii[] $A^{+}(\R)\oplus \R$, where $A^{+}(\R)$
is the group of orientation preserving affine maps on $\R$, 
\iii[]
$SOLV^{+},SOLV^{-}$ are Lie groups whose Lie algebra is solvable
and has 2-dim square, 
\iii[] $SE(2)$ and $SH(2)$ are the groups of orientation preserving
motions of Euclidean and Hyperbolic plane respectively, 
\iii[]
$SL(2)$ and $SU(2)$ are the three dimensional simple Lie groups.
\ei 

Moreover it is easy to show that in each of these cases but one
all left-invariant bracket generating distributions are equivalent by
automorphisms of the Lie algebra. The only case where there
exists two non-equivalent distributions is the Lie algebra
$\mathfrak{sl}(2)$. More precisely a 2-dimensional subspace of
$\mathfrak{sl}(2)$ is called \emph{elliptic} (\emph{hyperbolic})
if the restriction of the Killing form on this subspace is
sign-definite (sign-indefinite). Accordingly, we use notation
$SL_{e}(2)$ and $SL_{h}(2)$ to specify on which subspace the
sub-Riemannian structure on $SL(2)$ is defined.

For a left-invariant structure on a Lie group the invariants $\chi$ and $\kappa$ are constant functions and allow us to distinguish non isometric structures.
To complete the classification we can restrict ourselves to \emph{normalized} sub-Riemannian structures, i.e. structures that satisfy
\bqn\label{eqeq}\chi=\kappa=0, ~~~~\qquad \tx{or}\qquad~~~~ \chi^2+\kappa^2=1.\eqn
Indeed $\chi$ and $\kappa$ are homogeneous with respect to dilations of the orthonormal frame, that means rescaling of distances on the manifold. Thus we can always rescale our structure in such a way that \eqref{eqeq} is satisfied.

%This condition is not restrictive at all in our classification since, when we dilate the orthonormal frame, all distances on our manifold are multiplied by a constant and the new structure is a rescaling of the previous one.

To find missing discrete invariants, i.e. to distinguish between normalized structures with same $\chi$ and $\kappa$, we then show that it is always possible to select a canonical orthonormal frame for the sub-Riemannian structure such that all structure constants of the Lie algebra of this frame are invariant with respect to local isometries. Then the commutator relations of the Lie algebra generated by the canonical frame determine in a unique way the sub-Riemannian structure.

\begin{figure}[!htp]
\centering
\includegraphics[scale=0.62]{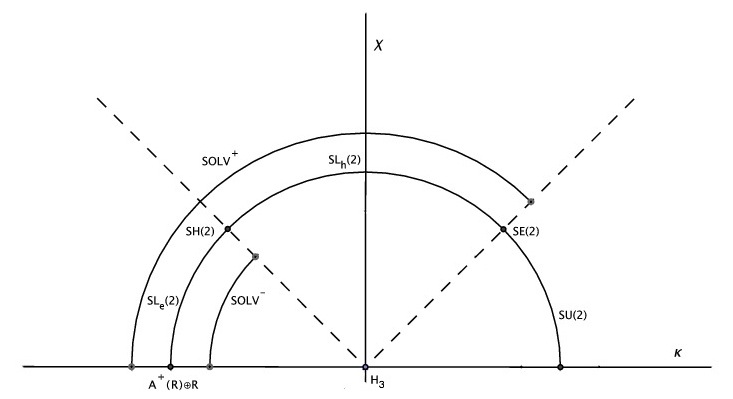}
\caption{Classification} \label{fig:classifica}
\end{figure}

Collecting together these results we prove the following
\bt \label{t-class} All left-invariant sub-Riemannian structures on 3D Lie groups are classified up to local isometries and dilations as in Figure \ref{fig:classifica}, where a structure is identified by the point $(\kappa,\chi)$ and two distinct points represent non locally isometric structures.

Moreover
\bi
\iii[$(i)$] If $\chi=\kappa=0$ then the structure is locally isometric to the Heisenberg group,
\iii[$(ii)$] If $\chi^{2}+\kappa^{2}=1$ then there exist no more than three non isometric normalized sub-Riemannian structures with these invariants; in particular there exists a unique normalized structure on a unimodular Lie group (for every choice of $\chi,\kappa$).
\iii[$(iii)$] If $\chi\neq0$ or $\chi=0,\kappa\geq0$, then two structures are locally isometric if and only if their Lie algebras are isomorphic.
\ei
\et

In other words every left-invariant sub-Riemannian structure is locally isometric to a normalized one that appear in Figure \ref{fig:classifica}, where we draw points on different circles 
since we consider equivalence classes of structures up to dilations. In this way it is easier 
to understand how many normalized structures there exist for some fixed value of the local invariants. Notice that unimodular Lie groups are those that appear in the middle circle (except for $A^{+}(\R)\oplus \R$).

From the proof of Theorem \ref{t-class} we get also a
uniformization-like theorem for \virg{constant curvature}
manifolds in the sub-Riemannian setting: 
\bc \label{t:uniforme}
Let $M$ be a complete simply connected 3D contact sub-Riemannian
manifold. Assume that $\chi=0$ and $\kappa$ is costant on $M$.
Then $M$ is isometric to a left-invariant sub-Riemannian
structure. More precisely: \bi \iii[(i)] if $\kappa=0$ it is
isometric to the Heisenberg group $H_3$, \iii[(ii)] if $\kappa=1$
it is isometric to the group $SU(2)$ with Killing metric,
\iii[(iii)] if $\kappa=-1$ it is isometric to the group
$\widetilde{SL}(2)$ with elliptic type Killing metric, \ei where
$\widetilde{SL}(2)$ is the universal covering of $SL(2)$. \ec

%We stress that in Corollary \ref{t:uniforme} we do not require to be left-invariant

Another byproduct of the classification is the fact that there exist non isomorphic Lie groups with locally isometric sub-Riemannian structures. Indeed, as a consequence of Theorem \ref{t-class}, we get that there exists a unique normalized left-invariant structure defined on $A^+(\R)\oplus\R$ having $\chi=0,\kappa=-1$. Thus $A^+(\R)\oplus\R$ is locally isometric to the group $SL(2)$ with elliptic type Killing metric by Corollary \ref{t:uniforme}.

This fact was already noted in \cite{falbel} as a consequence of the classification. In this paper we explicitly compute the global sub-Riemannian isometry between $A^+(\R)\oplus\R$ and the universal covering of $SL(2)$ by means of Nagano principle. We then show that this map is well defined on the quotient, giving a global isometry between the group $A^+(\R)\times S^1$ and the group $SL(2)$, endowed with the sub-Riemannian structure defined by the restriction of the Killing form on the elliptic distribution. 

The group $A^{+}(\R)\oplus \R$ can be interpreted as the subgroup of the affine maps on the plane %$\tx{Aff}(\R^{2})$ 
that acts as an orientation preserving affinity on one axis and as translations on the other one\footnote{We can recover the action as an affine map identifying $(x,y)\in \R^{2}$ with  $(x, y, 1)^{T}$ and
$$\begin{pmatrix}
a & 0 & b\\
0 & 1 & c\\
0  & 0  & 1\\
\end{pmatrix}\begin{pmatrix}
x\\ y\\ 1\\
\end{pmatrix}=\begin{pmatrix}
ax+b\\ y+c\\ 1\\
\end{pmatrix}.$$}

\begin{equation*}\label{eq:grupL2+R}
A^+(\mb{R})\oplus \mb{R}:=\left \{
\begin{pmatrix}
a & 0 & b\\
0 & 1 & c\\
0  & 0  & 1\\
\end{pmatrix}, \ a>0, b,c\in \mb{R}
    \right \}.
\end{equation*}

The standard left-invariant sub-Riemannian structure on $A^{+}(\R)\oplus \R$ is defined by the orthonormal frame $\distr=\tx{span}\{e_{2},e_{1}+e_{3}\}$, where 
\begin{equation*}
e_1=\begin{pmatrix}
0 & 0 & 1 \\
0 & 0 & 0 \\
0 & 0 & 0 \\
\end{pmatrix}, \quad
e_2=\begin{pmatrix}
-1 & 0 & 0\\
0 & 0 & 0\\
0 & 0 & 0 \\
\end{pmatrix}, \quad
e_3=\begin{pmatrix}
0 & 0 & 0\\
0 & 0 & 1\\
0 & 0 & 0 \\
\end{pmatrix},
\end{equation*}
is a basis of the Lie algebra of the group, satisfying $[e_{1},e_{2}]=e_{1}$.
%The group of orientation preserving affine maps on the real line can be described as the matrix subgroup
%\begin{align*}
%A^{+}(\R)=\left \{
%\begin{pmatrix}
%a & b \\
%0 & 1 \\
%\end{pmatrix}; \ a,b \in \R, a>0
%    \right \}
%\end{align*}
%and we can recover the action of an element as an affine map identifying $x\in \R$ with  $(x \ 1)^{T}$
%$$\begin{pmatrix}
%a & b \\
%0 & 1 \\
%\end{pmatrix}\begin{pmatrix}
%x  \\
%1 \\
%\end{pmatrix}=\begin{pmatrix}
%ax + b \\
%1 \\
%\end{pmatrix}$$

The subgroup $A^{+}(\R)$ is topologically homeomorphic to the half-plane
$\{(a,b)\in \R^{2}, a>0\}$ which can be descirbed in standard polar
coordinates as $\{(\rho,\theta)|\, \rho>0,-\pi/2 <
\theta< \pi/2 \}$.

%\ffoot{In \cite{laplacian} is studied laplacian
%on unimodular lie group and \ldots here we found a solution for
%nonunimodular\ldots}

\bt \label{t-sriso} The diffeomorphism $\Psi: A^+(\R)\times S^1 \longrightarrow SL(2)$ defined by
 \begin{gather}\label{eq:iso}
\Psi(\rho,\theta,\ph)=\frac{1}{\sqrt{\rho \cos\theta}}\begin{pmatrix}
\displaystyle{\cos \ph} &
\displaystyle{\sin \ph}\\
\displaystyle{\rho \sin (\theta - \ph)}&
\displaystyle{\rho \cos (\theta - \ph)}\\
\end{pmatrix},
\end{gather} where $(\rho, \theta)\in A^+(\R)$ and $\ph \in S^{1}$,
is a global sub-Riemannian isometry.
\et
Using this global sub-Riemannian isometry as a change of coordinates one can recover the geometry of the sub-Riemannian structure on the group $A^+(\R)\times S^1$, starting from the analogous properties of $SL(2)$ (e.g. explicit expression of the sub-Riemannian distance, the cut locus). In particular we notice that, since $A^+(\R)\times S^1$ is not unimodular, the canonical sub-Laplacian on this group is not expressed as a sum of squares. Indeed if $X_{1},X_{2}$ denotes the left-invariant vector fields associated to the orthonormal frame, the sub-Laplacian is expressed as follows
$$L_{sR}=X_{1}^{2}+X^{2}_{2}+X_{1}.$$
Moreover in the non-unimodular case the generalized Fourier transform method, used in \cite{laplacian}, cannot apply . Hence the heat kernel of the corresponding heat equation cannot be computed directly. On the other hand one can use the map \eqref{eq:iso} to express the solution in terms of the heat kernel on $SL(2)$.

\medskip
{\bf Acknowledgements.} The authors are grateful to the anonymous reviewer who noticed some important papers in reference to our work. 

\section{Basic definitions}\label{sec:bd}
We start recalling the definition of sub-Riemannian manifold.
\bdeff
A \emph{sub-Riemannian manifold} is a triple $(M,\distr,\metr)$,
where
\bi
\iii[$(i)$] $M$ is a smooth connected $n$-dimensional manifold,
\iii[$(ii)$] $\distr$ is a smooth distribution of constant rank $k< n$,
%satisfying the \emph{H\"ormander condition},
i.e. a smooth map that associates to every $q\in M$  a $k$-dimensional subspace $\distr_{q}$ of $T_qM$,
%and we have
%\bqn \label{Hor}
%\text{span}\{[X_1,[\ldots[X_{j-1},X_j]\ldots]](q)~|~X_i\in\mathrm{Vec}_H(M),\, j\in \N\}=T_qM, \quad \all q\in M
%\eqn
%where $\mathrm{Vec}_{H}(M)$ denotes the set of \emph{horizontal smooth vector fields} on $M$, i.e.
\iii[$(iii)$] $\metr_q$ is a Riemannian metric on
$\distr_{q}$, that is smooth with respect to $q\in M$. \ei \edeff
%
%Let $M$ be a smooth manifold. A \emph{sub-Riemannian structure of rank k} on $M$ is $(\distr,\la\cdot,\cdot\ra)$ where  $\distr\subset TM$ is a rank $k$ smooth distribution ($1< k< n$) and $\la\cdot,\cdot\ra$ is a smooth family of euclidean products $\la\cdot,\cdot\ra_q$ defined on $\distr_q$.
The set of smooth sections of the distribution
$$\overline{\distr}:=\{f\in \VecM | \, f(q)\in \distr_q,\ \all q\in M \}\subset \VecM.$$ is a subspace of the space  of the smooth vector fields on $M$ and its elements are said \emph{horizontal} vector fields.

A Lipschitz continuous curve $\g:[0,T]\to M$ is \emph{admissible} (or \emph{horizontal}) if its derivative is a.e. horizontal, i.e. if $\dot{\g}(t)\in \distr_{\g(t)}$ for a.e. $t\in [
0,T]$. We denote with $\Omega_{pq}$ the set of admissible paths joining $p$ to $q$.

Given an admissible curve $\g$ it is possible to define its \emph{lenght}
$$\ell(\g)=\int_0^T \sqrt{ \metr_{\g(t)}(\dot{\g}(t),\dot{\g}(t))}\,dt.$$

The \emph{Carnot-Caratheodory distance} induced by the sub-Riemannian structure is
$$d(p,q)=\inf \{\ell(\g), \g \in \Omega_{pq}\}.$$

In the following we always assume that the distribution $\distr$ satisfies the \emph{bracket generating condition} (also known as \emph{H$\ddot{o}$rmander condition}), i.e. the Lie algebra generated by the horizontal vector fields spans at every point the tangent space to the manifold
$$\tx{span}\{[f_1,\ldots,[f_{j-1},f_j]](q),f_i\in \overline{\distr}, j\in \N\}=T_qM,\quad \all q\in M.$$

Under this hypothesis the classical
Chow-Rashevsky Theorem \cite{chow,rashevsky} implies that $d$ is a well
defined metric on $M$ and it induces on $M$ the original topology.

\bdeff A \emph{sub-Riemannian isometry} between two sub-Riemannian manifolds $(M,\distr,\metr)$ and $(N,\distr^{\prime},\metr^{\prime})$ is a diffeomorphism $\phi:M\rightarrow N$ that satisfies
\begin{itemize}
\item[$(i)$] $\phi_*(\distr)=\distr^{\prime},$
 \item[$(ii)$] $\metr(
f_1,f_2) =\metr^{\prime}(\phi_*f_1,\phi_*f_2), \quad \all
\, f_1,f_2 \in \overline{\distr}.$
\end{itemize}
\edeff

\bdeff Let $M$ be a $2m+1$ dimensional manifold. A sub-Riemannian
structure on $M$ is said to be \emph{contact} if $\distr$ is a
contact distribution, i.e. $\distr= \ker \omega$, where $\omega\in\Lambda^{1}M$ satisfies $\left(\bigwedge^md\omega\right)\wedge\omega \neq 0$. Notice that a contact structure is forced to
be bracket generating.\edeff
The contact structure endows $M$ with a canonical orientation. On the other hand we will not fix an orientation on the distribution $\distr$.
%\ffoot{orientable?} and the contact 1-form, which is defined up to a smooth function, could be selected in such a way that $(d\omega)^n$ at poit $q$ coincide with the volume form on the euclidean plane $\distr_q$.

Now we briefly recall some facts about sub-Riemannian geodesics. In particular we define the sub-Riemannian Hamiltonian.
%\subsection{Geodesics}

Let $M$ be a sub-Riemannian manifold and fix $q_0\in M$. We define the \emph{endpoint map} (at time 1) as $$F:\mathcal{U}\to M, \quad F(\g)=\g(1),$$ where  $\mathcal{U}$ denotes  the set of admissible trajectories starting from $q_0$ and defined at time $t=1$. If we fix a point $q_1\in M$, the problem of finding shortest paths from $q_0$ to $q_1$ is equivalent to the following one
\bqn \label{eq:min}
\min_{F^{-1}(q_1)}J(\g),\qquad J(\g):=\frac{1}{2}\int_0^1|\dot{\g}(t)|^2dt,
\eqn
 where $J$ is the action functional. Indeed,
it is a standard fact that Cauchy-Schwartz inequality implies that an admissible curve realizes this minimum if and only if it is an arc-lenght parametrized $\ell$-minimizer.

Then the Lagrange multipliers rule implies that any solution of \eqref{eq:min} is either a critical point of $F$ or a solution of the equation
\bqn \label{eq:norm}
\lam_{1} D_{\g}F=d_{\g}J,\qquad \g\in \mathcal{U},
\eqn
 for some $\lam_{1} \in T^*_{\g(1)}M$.
Solutions of equation \eqref{eq:norm} are said \emph{normal geodesics} while critical ponits of $F$ are said \emph{abnormal geodesics}.

Now we can define the \emph{sub-Riemannian Hamiltonian} $h\in
C^{\infty}(T^*M)$ as follows: \bqn \label{eq:srham}
h(\lam)=\max_{u\in \distr_q} \{\la \lam,u\ra-\frac{1}{2}|u|^2\},
\qquad \lam\in T^{*}M, \quad q=\pi(\lam), \eqn where $\la\cdot,\cdot\ra$ denotes the
standard pairing between vectors and covectors. The Pontryagin
Maximum Principle gives a perfect characterization of our
geodesics. In fact it can be shown that in the contact case there are
no abnormal geodesics and a pair $(\g,\lam_1)$ satisfies
\eqref{eq:norm} if and only if there exists a curve $\lam(t) \in
T^*_{\g(t)}M$ that is a  solution of the Hamiltonian system
$\dot{\lam}(t)=\vec{h}(\lam(t))$ with boundary condition
$\lam(1)=\lam_1$.
%the optimal control problem (where $f_1,\ldots,f_k$ is a local orthonormal frame)
%\begin{align}
%\dot{\g}(t)&=\sum_{i=1}^k u_i(t) f_i(\g(t))  \notag\\
%\g(0)&=q_0, \quad \g(T)=q_1,\\
%\ell(\g)&\to \min \notag
%\end{align}

\brem Locally the sub-Riemannian structure can be given assigning a set of $k$ smooth linearly independent vector fields that are orthonormal
\bqn \label{eq:local}
\distr_q=\tx{span}\{f_1(q),\ldots,f_k(q)\},\qquad \metr_{q}(f_i(q),f_j(q))=\delta_{ij}.
\eqn
Notice that if we consider a new orthonormal frame which is a rotation of the previous one, we define the same sub-Riemannian structure.

Following this notation a \emph{local isometry} between two structures defined by the orthonormal frames
$\distr_M=\tx{span}(f_1,\ldots,f_k)$, $
 \distr_N=\tx{span}(g_1,\ldots,g_k)$ is given by a local diffeomorphism such that
$$\phi:M
\rightarrow N, \quad \phi_*(f_i)=g_i,\qquad \all\, i=1,\ldots,k.$$

In this setting admissible trajectories are solutions of the equation
$$\dot{\g}(t)=\sum_{i=1}^k u_i(t) f_i(\g(t)), \quad \text{for a.e. }t\in [0,T] ,$$
for some measurable and bounded control functions
$u(t)=(u_1(t),\ldots,u_k(t)), u_i \in L^{\infty}([0,T]).$
Lenght and action of this curve are expressed as follows
$$\ell(\g)=\int_0^T |u(t)|dt, \quad J(\g)=\frac{1}{2}\int_0^T|u(t)|^2dt,$$
where $|\cdot |$ denotes standard Euclidean norm in $\R^{k}$.

Moreover the sub-Riemannian Hamiltonian \eqref{eq:srham} is written as
$$h(\lam)=\frac{1}{2}\sum_{i=1}^{k} h_i^2(\lam), \qquad \text{ where }\quad  h_i(\lam)=\la \lam, f_i(q)\ra, \quad q=\pi(\lam).$$
Notice that $h_i:T^*M \to \R$ are smooth functions on $T^{*}M$ that are linear on fibers, associated to the vector fields of the frame. The sub-Riemannian Hamiltonian $h$ is a smooth function on $T^*M$ which contains all the informations about the sub-Riemannian structure. Indeed it does not depend on the orthonormal frame selected $\{f_{1},\ldots,f_{k}\}$, i.e. is invariant for rotations of the frame, and the annichilator of the distribution at a point $\distr_q^{\perp}$
%(more precisely, its annichilator $\distr^\perp_q$}
can be recovered as the kernel of the restriction of $h$ to the fiber $T^*_qM$
$$\ker h|_{T^*_qM}=\{\lam \in T^*_qM|\  h_i(\lam)=0, \ i=1,\ldots,k\}=\distr^\perp_q.$$
%\ffoot{$\distr^\perp_q:=\{\lam \in T^*_qM, \la\lam,v\ra=0, \forall v\in \distr_q\}$.}

 \erem

\brem\label{r:r}
A sub-Riemannian structure on a Lie group $G$ is said to be \emph{left-invariant} if $$\distr_{gh}=L_{g*}\distr_h, \qquad \la v,w \ra_h=\la L_{g*}v,L_{g*}w\ra_{gh}, \qquad \forall g,h\in G.$$
where $L_{g}$ denotes the left multiplication map on the group.
In particular, to define a left-invariant structure, it is sufficient to fix a subspace of the Lie algebra $\mathfrak{g}$ of the group and an inner product on it.

We also remark that in this case it is possible to have in \eqref{eq:local} a global equality, i.e. to select $k$ globally linearly independent orthonormal vector fields.
%\ffoot{In the following, when we write $\distr=\tx{span}(f_1,\ldots,f_k)$, we refer to the left-invariant structure where the vector fields $f_j$ are orthonormaland this is globally defined.}
\erem

\section{Sub-Riemannian invariants} \label{sec:inv}
In this section we study a contact sub-Riemannian structure on a 3D manifold and we give a brief description of its two invariants (see also \cite{agrexp}). We start with the following characterization of contact distributions.
\bl \label{lem:vol}Let $M$ be a 3D manifold, $\omega \in \Lambda^1M$ and $\distr=\ker \omega$. The following are equivalent:
\bi
\iii[$(i)$] $\distr$ is a contact distribution,
\iii[$(ii)$] $d\omega \big|_{\distr}\neq0$,
\iii[$(iii)$] $\all f_1,f_2\in \overline{\distr}$ linearly independent, then $[f_1,f_2]\notin \overline{\distr}$.
\ei
Moreover, in this case, the contact form can be selected in such a way that $d\omega \big|_{\distr}$ coincide with the Euclidean volume form on $\distr$.
\el
By Lemma \ref{lem:vol} it is not restrictive to assume that
the sub-Riemannian structure satisfies:
\begin{gather}
(M,\omega) \text{ is a 3D contact structure}, \notag  \\
\distr=\tx{span}\{f_1,f_2\}=\ker \omega,\label{eq:setting}\\
\metr(f_i,f_j)=\delta_{ij}, \quad
d\omega(f_1,f_2)=1. \notag
\end{gather}
We stress that in \eqref{eq:setting} the orthonormal frame
$f_1,f_2$ is not unique. Indeed every rotated frame (where the angle
of rotation depends smoothly on the point) defines the same
structure.

The sub-Riemannian Hamiltonian \eqref{eq:srham} is written
 \bqn \label{eq:ham3d}
 h=\frac{1}{2}(h_1^2+h_2^2). \notag
 \eqn

\bdeff In the setting \eqref{eq:setting} we define the \emph{Reeb vector field} associated to the contact structure as the unique vector field $f_0$ such that
\begin{align} \label{eq:deff0}
\omega(f_0)&=1, \notag\\
d\omega(f_0,\cdot)&=0.
\end{align}\edeff
From the definition it is clear that $f_0$ depends only on the sub-Riemannian structure (and its orientation) and not on the frame selected.

Condition \eqref{eq:deff0} is equivalent to%\ffoot{$d\omega(X,Y)=X\omega(Y)-Y\omega(X)-\omega([X,Y])\qquad \forall \, X,Y \in \VecM, \omega\in \Lambda^1(M), $}
\begin{align*}%\label{eq:deff0comm}
&[f_1,f_0],[f_2,f_0]\in \overline{\distr}, \notag\\
&[f_2,f_1]=f_0 \quad (\text{mod } \overline{\distr}).
\end{align*}
and we deduce the following expression for the Lie algebra of vector fields generated by $f_0,f_1,f_2$
%\footnote{$c_{ij}^k$ are said \emph{structure constant}}
\begin{align} \label{eq:algebracampi}
[f_1,f_0]&=c_{01}^1 f_1+c_{01}^2 f_2, \notag\\
[f_2,f_0]&=c_{02}^1 f_1+c_{02}^2 f_2, \\
[f_2,f_1]&=c_{12}^1 f_1+c_{12}^2 f_2+f_0, \notag
\end{align}
where $c_{ij}^{k}$ are functions on the manifold, called structure constants of the Lie algebra.

If we denote with $(\nu_0,\nu_1,\nu_2)$ the basis of $1$-form dual to $(f_{0},f_{1},f_{2})$, we can rewrite \eqref{eq:algebracampi} as:
\begin{align}\label{eq:1-formeduali}
d\nu_0&= ~\quad \nu_1\wedge \nu_2, \notag\\
d\nu_1&= c_{01}^1 \nu_0 \wedge \nu_1 +c_{02}^1\nu_0\wedge \nu_2 +c_{12}^1 \nu_1\wedge \nu_2, \\
d\nu_2&= c_{01}^2 \nu_0 \wedge \nu_1 +c_{02}^2\nu_0\wedge
\nu_2+c_{12}^2 \nu_1 \wedge \nu_2, \notag
\end{align}

Let $h_0(\lam)=\langle \lam,f_0(q) \rangle$ denote the  Hamiltonian linear on fibers
associated with the Reeb field $f_0$. We now compute
%, being $h$ and $h_{0}$ defined only by the sub-Riemannian structure,
the Poisson bracket $\{h,h_0\}$, denoting with $\{h,h_0\}_q$ its
restriction to the fiber $T^*_qM$. \bp \label{pr:poiss} The
Poisson bracket  $\{h,h_0\}_q$ is a quadratic form. Moreover we
have
\begin{gather}
\{h,h_0\}=c_{01}^1h_1^2+(c_{01}^2+c_{02}^1)h_1h_2+c_{02}^2h_2^2, \label{eq:parentesihh0}\\
c_{01}^1+c_{02}^2=0. \label{eq:traccia0}
\end{gather}
In particular, $\Delta^\perp_q\subset\ker\,\{h,h_0\}_q$ and
$\{h,h_0\}_q$ is actually a quadratic form on
$T^*_qM/\Delta^\perp_q=\Delta^*_q$. \ep

\begin{proof}
Using the equality $\{h_i,h_j\}(\lam)=\la\lam,[f_i,f_j](q)\ra$  we get
\begin{align*}
\{h,h_0\}&=\frac{1}{2}\{h_1^2+h_2^2,h_0\}=h_1\{h_1,h_0\}+h_2\{h_2,h_0\}\\
         %&=\langle p,f_1(q) \rangle \langle p,[f_1,f_0](q)\rangle
         %+\langle p,f_2(q) \rangle \langle p,[f_2,f_0](q)\rangle\\
         &=h_1(c_{01}^1h_1+c_{01}^2h_2)+h_2(c_{02}^1h_1+c_{02}^2h_2)\\
         &=c_{01}^1h_1^2+(c_{01}^2+c_{02}^1)h_1h_2+c_{02}^2h_2^2.
\end{align*}
Differentiating the first equation in
\eqref{eq:1-formeduali} we find:
\begin{align*}
0=d^2\nu_0&=d\nu_1 \wedge \nu_2 - \nu_1 \wedge d\nu_2 \\
        &=(c_{01}^1+c_{02}^2)\nu_0 \wedge \nu_1 \wedge \nu_2.
\end{align*}
which proves \eqref{eq:traccia0}.
\end{proof}

Being $\{h,h_0\}_q$ a quadratic form on the Euclidean plane $\distr_{q}$ (using the canonical identification of the vector space $\distr_{q}$ with its dual $\distr^{*}_{q}$ given by the scalar product), it is a standard fact that it can be interpreted as a symmetric operator on the plane itself. In particular its determinant and its trace are well defined. From \eqref{eq:traccia0} we get
$$\tx{trace}\, \{h,h_0\}_q=0.$$
It is natural then to define our \emph{first invariant} as the positive eigenvalue of this operator, namely:
\bqn \label{eq:defchi}
 \chi(q)=\sqrt{-\tx{det}\{h,h_0\}_q}.
\eqn
\brem \label{r:chipositivo} Notice that, by definition $\chi\geq0$, and it vanishes everywhere if and only if the flow of the Reeb vector field $f_{0}$ is a flow of sub-Riemannian isometries for $M$.
\erem
The \emph{second invariant}, which was found in \cite{agrexp} as a term of the asymptotic
expansion of conjugate locus, is defined in the following way
\bqn \label{eq:defkappa}
\kappa(q)=f_2(c_{12}^1)-f_1(c_{12}^2)-(c_{12}^1)^2-(c_{12}^2)^2+
\dfrac{c_{01}^2-c_{02}^1}{2}.
\eqn     where we refer to notation \eqref{eq:algebracampi}.
%We will refer also at $\kappa$ as \emph{sub-Riemannian curvature}, and we will see that it is not just a matter of notation.\ffoot{indeed
%there are some relation between $\kappa$ and Gaussian curvature (vedere se mettere la formula\ldots), see also Theorem \ref{t:uniforme}.}
A direct calculation shows that $\kappa$ is preserved by rotations
of the frame $f_1,f_2$ of the distribution, hence it depends only
on the sub-Riemannian structure.

$\chi$ and $\kappa$ are functions defined on the manifold; they
reflect intrinsic geometric properties of the sub-Riemannian
structure and are preserved by the sub-Riemannian isometries. In
particular, $\chi$ and $\kappa$ are constant functions for
left-invariant structures on Lie groups (since left translations
are isometries).

\section{Canonical Frames} \label{sec:frames}
In this section we want to show that it is always possible to
select a canonical orthonormal frame for the sub-Riemannian
structure. In this way we are able to find missing discrete
invariants and to classify sub-Riemannian structures simply
knowing structure constants $c_{ij}^k$ for the canonical frame.
We study separately the two cases $\chi \neq 0$ and $\chi=0$.

We start by rewriting and improving Proposition \ref{pr:poiss}
when $\chi \neq 0$. \bp \label{p:framenon0} Let $M$ be a 3D
contact sub-Riemannian manifold and $q\in M$. If $\chi(q)\neq 0$,
then there exists a local frame such that %\eqref{eq:parentesihh0} becomes: 
\bqn \label{eq:parconchi} \{h,h_0\}=2\chi h_1h_2. \eqn In
particular, in the Lie group case with left-invariant stucture,
there exists a unique (up to a sign) canonical frame
$(f_0,f_1,f_2)$ such that %\eqref{eq:algebracampi} become
\begin{align}\label{eq:algchinon0}
    [f_1,f_0]&=  c_{01}^2f_2 ,  \notag\\
    [f_2,f_0]&=c_{02}^1f_1, \\
    [f_2,f_1]&=c_{12}^1f_1  +c_{12}^2f_2+f_0. \notag
\end{align}
Moreover we have
\bqn \label{eq:invnellec}
\chi=\frac{c_{01}^2+c_{02}^1}{2}, \quad~~~~~ \kappa=-(c_{12}^1)^2-(c_{12}^2)^2+
\frac{c_{01}^2-c_{02}^1}{2}.
\eqn
\ep
\begin{proof} From Proposition \ref{pr:poiss} we know that the Poisson bracket
$\{h,h_0\}_q$ is a non degenerate symmetric operator with zero
trace. Hence we have a well defined, up to a sign, orthonormal frame by setting 
$f_1,f_2$ as the orthonormal isotropic vectors of this operator
(remember that $f_0$ depends only on the structure and not on the
orthonormal frame on the distribution). It is easily seen that in
both of these cases we obtain the expression \eqref{eq:parconchi}.
\end{proof}
\brem\label{r:13}
Notice that, if we change sign to $f_{1}$ or $f_{2}$, then $c_{12}^{2}$ or $c_{12}^{1}$, respectively, change sign in \eqref{eq:algchinon0}, while $c_{02}^{1}$ and $c_{01}^{2}$ are unaffected. Hence equalities \eqref{eq:invnellec} do not depend on the orientation of the sub-Riemannian structure. 
%since  n the case of left-invariant structures, since all structure constant in are constant functions, we could select, for instance, the frame in which  are positive and follows.
\erem

If $\chi=0$ the above procedure cannot apply. Indeed both trace and determinant of the operator vanish, hence we have $\{h,h_{0}\}_{q}=0$.
From \eqref{eq:parentesihh0} we get the identities
\begin{equation}\label{eq:chi=0}
c_{01}^1=c_{02}^2=0, \quad
    c_{01}^2+c_{02}^1=0.
\end{equation}
so that commutators \eqref{eq:algebracampi} simplify in (where $c=c_{01}^2$)
\begin{align}\label{eq:commutatorichi0}
    [f_1,f_0]&=  cf_2,   \notag\\
    [f_2,f_0]&=-cf_1, \\
    [f_2,f_1]&=c_{12}^1f_1  +c_{12}^2f_2  +f_0. \notag
\end{align}

We want to show, with an explicit construction, that also in this
case there always exists a rotation of our frame, by an angle that smoothly
depends on the point, such that in the new frame $\kappa$
is the only structure constant which appear in
\eqref{eq:commutatorichi0}.

%We begin with a useful lemma
\begin{lemma}\label{lem:rotazione}
Let $f_1,f_2$ be an orthonormal frame on $M$. If we denote with
$\widehat{f}_1,\widehat{f}_2$ the frame obtained from the previous
one with a rotation by an angle $\theta(q)$ and with
$\widehat{c}_{ij}^k$ structure constants of rotated frame, we
have:
\begin{align*}%\label{eq:rotaz}
\widehat{c}_{12}^1= \cos \theta (c_{12}^1-f_1(\theta))- \sin
\theta (c_{12}^2-f_2(\theta)), \\
\widehat{c}_{12}^2= \sin \theta (c_{12}^1-f_1(\theta))+ \cos
\theta (c_{12}^2-f_2(\theta)).
\end{align*}
\end{lemma}

Now we can prove the main result of this section.
\bp \label{t:rotaz} Let $M$ be a 3D simply connected contact sub-Riemannian manifold
such that $\chi=0$. Then there exists a rotation of the original frame $\widehat{f}_1,\widehat{f}_2$ such that:
\begin{align}\label{eq:canonicochi0}
    [\widehat{f}_1,f_0]&=\kappa \widehat{f}_2,   \notag\\
    [\widehat{f}_2,f_0]&=-\kappa\widehat{f}_1, \\
    [\widehat{f}_2,\widehat{f}_1]&= f_0. \notag
\end{align}
\ep
\begin{proof}  Using Lemma \ref{lem:rotazione} we can rewrite the statement in the
following way: there exists a function
$\theta:M\rightarrow \mb{R}$ such that
\begin{equation}\label{eq:fugualec}
f_1(\theta)=c_{12}^1,\quad f_2(\theta)=c_{12}^2.
\end{equation}
Indeed, this would imply $\widehat{c}_{12}^1=\widehat{c}_{12}^2=0$ and $\kappa=c$.

Let us introduce simplified notations $c_{12}^1=\alpha_1,\
c_{12}^2=\alpha_2$. Then \bqn \label{eq:1}
\kappa=f_2(\alpha_1)-f_1(\alpha_2)-(\alpha_1)^2-(\alpha_2)^2+c.
\eqn If  $(\nu_0,\nu_1,\nu_2)$ denotes the dual basis to
$(f_0,f_1,f_2)$ we have
$$d\theta= f_0(\theta)\nu_0+f_1(\theta)\nu_1+f_2(\theta)\nu_2.$$
and from \eqref{eq:commutatorichi0} we get:
\begin{align*}
f_0(\theta)&=([f_2,f_1]-\alpha_1f_1-\alpha_2f_2)(\theta) \\
           &=f_2(\alpha_1)-f_1(\alpha_2)-\alpha_1^2-\alpha_2^2 \\
           &=\kappa-c.
\end{align*}
Suppose now that \eqref{eq:fugualec} are satisfied, we get
\begin{equation} \label{eq:unoforma}
d\theta=(\kappa-c)\nu_0+\alpha_1\nu_1+\alpha_2\nu_2=:\eta.
\end{equation}
with the r.h.s. independent from $\theta$.

To prove the theorem we have to show that $\eta$ is an exact 1-form. Since the manifold is simply connected, it is sufficient to prove that $\eta$ is closed.
If we denote $\nu_{ij}:=\nu_i \wedge \nu_j$ dual equations of \eqref{eq:commutatorichi0} are:
\begin{align*}
d\nu_0&=\nu_{12},\\
d\nu_1&=-c\nu_{02} +\alpha_1 \nu_{12}, \\
d\nu_2&=c \nu_{01} -\alpha_2 \nu_{12}.
\end{align*}
and differentiating we get two nontrivial relations:
\begin{gather}\label{eq:struttchi01}
f_1(c)+c\alpha_2+f_0(\alpha_1)=0,\\
 f_2(c)-c\alpha_1+f_0(\alpha_2)=0.\label{eq:struttchi02}
\end{gather}
Recollecting all these computations we prove the closure of $\eta$
\begin{align*}
d\eta &= d(\kappa-c)\wedge \nu_0+(\kappa-c)d\nu_0+d\alpha_1 \wedge \nu_1+\alpha_1d\nu_1+d\alpha_2 \wedge \nu_2+\alpha_2d\nu_2\\
        &= -dc \wedge \nu_0+(\kappa-c)\nu_{12}+\\
        &\qquad +f_0(\alpha_1)\nu_{01}-f_2(\alpha_1)\nu_{12}+
           \alpha_1(\alpha_1 \nu_{12}-c\nu_{02})\\
        &\qquad \qquad +f_0(\alpha_2)\nu_{02}+f_1(\alpha_2)\nu_{12} +\alpha_2(c \nu_{01} -\alpha_2 \nu_{12})   \\
        &= (f_0(\alpha_1)+ \alpha_2 c +f_1(c)) \nu_{01} \\
        &\qquad +(f_0(\alpha_2)- \alpha_1 c +f_2(c)) \nu_{02}\\
        &\qquad \qquad
        +(\kappa-c-f_2(\alpha_1)+f_1(\alpha_2)+\alpha_1^2+\alpha_2^2)\nu_{12}\\
        &=0.
\end{align*}
where in the last equality we use \eqref{eq:1} and \eqref{eq:struttchi01}-\eqref{eq:struttchi02}.
\end{proof}

\section{Proof of Theorem \ref{t-class}} \label{sec:class}
Now we use the results of the previous sections to prove Theorem \ref{t-class}.

In this section $G$ denotes a 3D Lie group, with Lie algebra $\mathfrak{g}$, endowed with a left-invariant sub-Riemannian structure defined by the orthonormal frame $f_{1},f_{2}$, i.e. $$\distr= \text{span}\{f_1,f_2\}\subset \mathfrak{g}, \qquad \text{span}\{f_1,f_2,[f_1,f_2]\}=\mathfrak{g}.$$
Recall that for a 3D left-invariant structure to be bracket generating is equivalent to be contact, moreover the Reeb field $f_0$ is also a left-invariant vector field by construction.

%In the left-invariant case, the local geometry of the sub-Riemannian structure on $G$ is determined by the Lie algebra $\mathfrak{g}$.

From the fact that, for left-invariant structures, local
invariants are constant functions (see Remark \ref{r:r}) we
obtain a necessary condition for two structures to be locally isometric.

\bp Let $G,H$ be 3D Lie groups with locally isometric
sub-Riemannian structures. Then $\chi_G=\chi_H$ and
$\kappa_G=\kappa_H$. \ep Notice that this condition is not
sufficient. It turns out that there can be up to three mutually
non locally isometric normalized structures with the same invariants $\chi,
\kappa$.
% and in Section \ref{sec:iso} we will explicitly construct a remarkable example of this fact.

\brem  It is easy to see that $\chi$ and $\kappa$ are homogeneous of degree 2 with respect to dilations of the frame. Indeed assume that the sub-Riemannian structure $(M,\distr,\metr)$ is locally defined by the orthonormal frame $f_{1},f_{2}$, i.e.
$$\distr=\tx{span}\{f_1,f_2\},\qquad \metr(f_{i},f_{j})=\delta_{ij}.$$
Consider now the dilated structure $(M,\distr, \til{\metr})$ defined by the orthonormal frame $\lam f_{1},\lam f_{2}$
$$\distr=\tx{span}\{f_1,f_2\},\qquad \til{\metr}(f_{i},f_{j})=\frac{1}{\lam^{2}}\delta_{ij}, \qquad \lam>0.$$
If $\chi, \kappa$ and $\til{\chi}, \til{\kappa}$ denote the invariants of the two structures respectively, 
%In other words, if  we consider the two orthonormal frames on the same manifold
%$$\distr=\tx{span}\{f_1,f_2\},\qquad \widetilde{\distr}=\tx{span}\{\lambda f_1,\lambda f_2\}, \qquad \lam \in \R^{+}$$
%and denote $\til{\chi}$ and $\til{\kappa}$ local invariants for
%the dilated structure, 
we find$$\widetilde{\chi}=\lambda^2
\chi,\qquad  \widetilde{\kappa}=\lambda^2 \kappa, \qquad \lam>0.$$ 
A dilation of the orthonormal frame corresponds to a multiplication by a
factor $\lam>0$ of all distances in our manifold. Since we are interested in
a classification by local isometries, we can always suppose (for a
suitable dilation of the orthonormal frame) that the local invariants
of our structure satisfy
$$\chi=\kappa=0, \qquad \text{or} \qquad \chi^2+\kappa^2=1,$$
and we study equivalence classes with respect to local isometries.
\erem
Since $\chi$ is non negative by definition (see Remark \ref{r:chipositivo}), we study separately the two cases $\chi > 0$ and $\chi=0$.

\subsection{Case $\chi > 0$}
Let $G$ be a 3D Lie group with a left-invariant sub-Riemannian
structure such that $\chi \neq 0$. From Proposition
\ref{p:framenon0} we can assume that $\distr=
\text{span}\{f_1,f_2\}$ where $f_1,f_2$ is the canonical frame of
the structure. From  \eqref{eq:algchinon0} we obtain the dual
equations
\begin{align}\label{eq:1-formedualichinon0}
d\nu_0&= \nu_1\wedge \nu_2, \notag\\
d\nu_1&= c_{02}^1 \nu_0\wedge \nu_2 +c_{12}^1 \nu_1\wedge \nu_2, \\
d\nu_2&= c_{01}^2 \nu_0 \wedge \nu_1 +c_{12}^1 \nu_1 \wedge
\nu_2.\notag
\end{align}
Using $d^2=0$ we obtain structure equations
%\begin{align*}
%0=d^2\nu_1&= c_{02}^1c_{12}^2\, \nu_0 \wedge \nu_1 \wedge \nu_2\\
%0=d^2\nu_2&= c_{01}^2c_{12}^1\, \nu_0 \wedge \nu_1 \wedge \nu_2
%\end{align*}
%which we can summarize in
\begin{equation}\label{eq:sistema}
\left\{%
\begin{array}{ll}
    c_{02}^1c_{12}^2=0, \\[0.2cm]
    c_{01}^2c_{12}^1=0.\\
\end{array}%
\right.
\end{equation}

%Now we are able to classify all left-invariant structures with $\chi \neq 0$.
We know that the structure constants of the canonical frame are invariant by local isometries (up to change signs of $c_{12}^{1},c_{12}^{2}$, see Remark \ref{r:13}). Hence, every different choice of coefficients in \eqref{eq:algchinon0} which satisfy also \eqref{eq:sistema}  will belong to a different class of non-isometric structures.

Taking into account that $\chi> 0$ implies that $c_{01}^2$ and $c_{02}^1$ cannot be both non positive (see \eqref{eq:invnellec}), we have the following cases:
\bi
\item[$(i)$] $c_{12}^1=0$ and $c_{12}^2=0$.
In this first case we get
\begin{align*}
    [f_1,f_0]&=c_{01}^2f_2, \\
    [f_2,f_0]&=c_{02}^1f_1, \\
    [f_2,f_1]&=f_0,
\end{align*}
and formulas \eqref{eq:invnellec} imply
$$\chi=\frac{c_{01}^2+c_{02}^1}{2}>0, \qquad \kappa=\frac{c_{01}^2-c_{02}^1}{2}.$$
In addition, we find the relations between the invariants
\begin{equation}\label{eq:relazchikappa}
    \chi+\kappa=c_{01}^2,\qquad \chi-\kappa=c_{02}^1.    \notag
\end{equation}
We have the following subcases: \bi \iii[$(a)$] If $c_{02}^1=0$ we
get the Lie algebra  $\mathfrak{se}(2)$ of the group $SE(2)$ of
the Euclidean isometries of $\R^2$, and it holds $\chi=\kappa$.
\iii[$(b)$] If $c_{01}^2=0$ we get the Lie algebra
$\mathfrak{sh}(2)$ of the group $SH(2)$ of the Hyperbolic
isometries of $\R^2$, and it holds $\chi=-\kappa$. \iii[$(c)$] If
$c_{01}^2>0$ and $ c_{02}^1<0$ we get the Lie algebra
$\mathfrak{su}(2)$ and $\chi-\kappa<0$. \iii[$(d)$] If
$c_{01}^2<0$ and $ c_{02}^1>0$ we get the Lie algebra
$\mathfrak{sl}(2)$ with $\chi+\kappa<0$. \iii[$(e)$] If
$c_{01}^2>0$ and $ c_{02}^1>0$ we get the Lie algebra
$\mathfrak{sl}(2)$ with $\chi+\kappa>0,\chi-\kappa>0$. \ei
\vspace{0.3cm}
\item[$(ii)$] $c_{02}^1=0$ and $c_{12}^1=0$.
In this case we have
\begin{align}\label{eq:solv+}
           [f_1,f_0]&=c_{01}^2f_2,  \notag \\
           [f_2,f_0]&=0, \\
           [f_2,f_1]&=c_{12}^2f_2  +f_0, \notag
\end{align}
and necessarily $c_{01}^2 \neq 0$. Moreover we get
$$\chi=\frac{c_{01}^2}{2}>0, \qquad \kappa=-(c_{12}^2)^2+\frac{c_{01}^2}{2},$$
from which it follows
$$\chi- \kappa\geq0.$$
The Lie algebra $\mathfrak{g}=\text{span}\{f_{1},f_{2},f_{3}\}$ defined by \eqref{eq:solv+} satisfies dim$\,[\mathfrak{g},\mathfrak{g}]=2$, hence it can be interpreted as the operator $A=\text{ad }f_1$ which acts on the subspace span$\{f_0,f_2\}$. Moreover, it can be easily computed that
$$\text{trace }A=-c_{12}^2,\qquad \det
A=c_{01}^2>0,$$
and we can find the useful relation
\begin{equation}\label{eq:relazcaso1}
2\frac{\tx{trace}^2A}{\det A}=1-\frac{\kappa}{\chi}.
\end{equation}

\item[$(iii)$] $c_{01}^2=0$ and $c_{12}^2=0$.
In this last case we get
\begin{align} \label{eq:solv-}
    [f_1,f_0]&=0, \notag\\
    [f_2,f_0]&=c_{02}^1f_1, \\
    [f_2,f_1]&=c_{12}^1f_1+f_0, \notag
\end{align}
and $c_{02}^1 \neq 0$. Moreover we get
$$\chi=\frac{c_{02}^1}{2}>0, \qquad \kappa=-(c_{12}^1)^2-\frac{c_{02}^1}{2},$$
from which it follows
$$\chi+ \kappa\leq0.$$
As before, the Lie algebra $\mathfrak{g}=\text{span}\{f_{1},f_{2},f_{3}\}$ defined by \eqref{eq:solv-} has two-dimensional square and it can be interpreted as the operator $A=\text{ad }f_2$ which acts on the plane span$\{f_0,f_1\}$. It can be easily seen that it holds
$$\text{trace }A=c_{12}^1,\qquad \det
A=-c_{02}^1<0,$$
and we have an analogous relation
\begin{equation}\label{eq:relazcaso2}
2\frac{\tx{trace}^2A}{\det A}=1+\frac{\kappa}{\chi}.
\end{equation}
\ei
\brem Lie algebras of cases $(ii)$ and $(iii)$ are \emph{solvable} algebras and we will denote respectively $\mathfrak{solv}^+$ and $\mathfrak{solv}^-$, where the sign depends on the determinant of the operator it represents.
In particular, formulas \eqref{eq:relazcaso1} and \eqref{eq:relazcaso2} permits to recover the ratio between invariants (hence to determine a unique normalized structure) only from intrinsic properties of the operator. Notice that if $c_{12}^2=0$ we recover the normalized structure $(i)$-$(a)$ while   if $c_{12}^1=0$ we get the case $(i)$-$(b)$.
\erem
\brem
The algebra $\mathfrak{sl}(2)$ is the only case where we can define two nonequivalent distributions which corresponds to the case that Killing form restricted on the distribution is positive definite (case $(d)$) or indefinite (case $(e)$). We will refer to the first one as the \emph{elliptic} structure on $\mathfrak{sl}(2)$, denoted
$\mathfrak{sl}_e(2)$, and with \emph{hyperbolic} structure in the other case, denoting $\mathfrak{sl}_h(2)$.
\erem
\subsection{Case $\chi =0$}

A direct consequence of Proposition \ref{t:rotaz} for left-invariant structures is the following
\bc \label{c:chi0kappadecide}  Let $G,H$ be Lie groups with left-invariant sub-Riemannian structures and assume $\chi_G=\chi_H=0$. Then $G$ and  $H$ are locally isometric if and only if $\kappa_G=\kappa_H$.
\ec
Thanks to this result it is very easy to complete our classification. Indeed it is sufficient to find all left-invariant structures such that $\chi=0$ and to compare their second invariant $\kappa$.

A straightforward calculation leads to the following list of the
left-invariant structures on simply connected three dimensional
Lie groups with $\chi=0$:
\begin{itemize}
\item[-] $H_{3}$ is the Heisenberg nilpotent group; then $\kappa=0$.
\item[-] $SU(2)$ with the Killing inner product; then
$\kappa>0$.
\item[-] $\widetilde{SL}(2)$ with the elliptic distribution and Killing
inner product; then $\kappa<0$.
\item[-] $A^+(\R)\oplus\R$; then $\kappa<0$.
\end{itemize}

\brem In particular, we have the following:
\begin{itemize}
\item[$(i)$] All left-invariant sub-Riemannian structures  on
$H_3$ are locally isometric,
\item[$(ii)$] There exists on
$A^+(\mb{R})\oplus \mb{R}$ a unique (modulo dilations) left-invariant sub-Riemannian structure, which is locally isometric to $SL_{e}(2)$ with the Killing metric.
\end{itemize}
\erem

Proof of Theorem \ref{t-class} is now completed and we can recollect our result as in Figure \ref{fig:classifica}, where we associate to every normalized structure a point in the $(\kappa,\chi)$ plane: either $\chi=\kappa=0$, or $(\kappa,\chi)$ belong to the semicircle
$$\{(\kappa,\chi)\in \R^2, \chi^2+\kappa^2=1, \chi>0\}.$$
Notice that different points means that sub-Riemannian structures are not locally isometric.

\section{Proof of Theorem \ref{t-sriso}} \label{sec:iso}
In this section we want to write explicitly the sub-Riemannian
isometry between $SL(2)$ and $A^+(\mb{R})\times S^1$.
%As a result of classification Theorem \ref{t-class} we get that groups $A^+(\R)\oplus \R$ and $SL(2)$ with the elliptic structure are locally isometric because they have the same sub-Riemannian curvature. Moreover the isometry could be extended to a global one considering the universal covering $\til{SL(2)}$ of $SL(2)$.
%It turns out that two nonisomorphic Lie groups have isomorphic sub-Riemannian structures, and this will give, passing to a quotient,
%a global isomorphism between $A^+(\R) \times S^1$ and $SL(2)$.

Consider the Lie algebra
$\mathfrak{sl}(2)=\{A\in M_2(\R),\  \tx{trace}(A)=0\}=\tx{span}\{g_1,g_2,g_3\}$,
where $$g_1=\frac{1}{2}\begin{pmatrix}
1 & 0 \\
0 & -1 \\
\end{pmatrix}, \quad g_2=\frac{1}{2}\begin{pmatrix}
0 & 1\\
1 & 0\\
\end{pmatrix},\quad  g_3=\frac{1}{2}\begin{pmatrix}
0 & 1\\
-1 & 0 \\
\end{pmatrix}.$$
The sub-Riemannian structure on $SL(2)$ defined by the Killing form on the elliptic distribution is given by the orthonormal frame
\bqn \label{eq:slsl}\Delta_{\mathfrak{sl}}=\tx{span}\{g_1,g_2\}, \qquad \text{and} \qquad g_0:=-g_3,\eqn is the Reeb vector field.
Notice that this frame is already canonical since equations \eqref{eq:canonicochi0} are satisfied. Indeed
$$[g_1,g_0]=-g_2=\kappa g_{2}.$$
Recall that the universal covering of $SL(2)$, which we denote
$\widetilde{SL}(2)$, is a simply connected Lie group with Lie algebra
$\mathfrak{sl}(2)$. Hence \eqref{eq:slsl} define a left-invariant
structure also on the universal covering.

On the other hand we consider the following coordinates on the Lie group $A^+(\mb{R})\oplus \mb{R}$, that are well-adapted for our further calculations 
 \begin{equation}\label{eq:grupL2+R}
A^+(\mb{R})\oplus \mb{R}:=\left \{
\begin{pmatrix}
-y & 0 & x\\
0 & 1 & z\\
0  & 0  & 1\\
\end{pmatrix}, \quad y<0, x,z\in \mb{R}
    \right \}.
\end{equation}
It is easy to see that, in these coordinates, the group law reads
\begin{equation}\label{eq:operazinA2coord}
(x,y,z)(x^{\prime},y^{\prime},z^{\prime})=(x-yx^{\prime},-yy^{\prime},z+z^{\prime})\notag,
\end{equation}
and its Lie algebra $\mathfrak{a}(\R)\oplus \R$ is generated by the vector fields
\begin{equation}\label{eq:campiA2}
e_1=-y\partial_x,\quad e_2=-y\partial_y,\quad e_3=\partial_z \notag,
\end{equation} with the only nontrivial commutator relation
$[e_1,e_2]=e_1$.

 The left-invariant structure on $A^+(\mb{R})\oplus \mb{R}$ is defined by the orthonormal frame \begin{align}\label{eq:coeff}
\distr_{\mathfrak{a}}&= \text{span}\{f_1,f_2\}, \notag \\
 f_1&:=e_2=-y\partial_y, \\
 f_2&:=e_1+e_3=-y\partial_x+\partial_z.\notag\end{align}
With straightforward calculations we compute the Reeb vector field
$
 f_0=-e_3=-\partial_z $.%\qquad [f_1,f_0]=-f_2$$

 This frame is not canonical since it does not satisfy equations \eqref{eq:canonicochi0}. Hence we can apply Proposition \ref{t:rotaz} to find the canonical frame, that will be no more left-invariant.

Following the notation of Proposition \ref{t:rotaz} we have \bl
The canonical orthonormal frame on $A^+(\mb{R})\oplus \mb{R}$ has
the form:
\begin{align}\label{eq:nuovoframeA2}
\widehat{f}_1&=\quad y \sin z \,
\partial_x-y \cos z \, \partial_y-\sin z \,\partial_z, \notag\\
\widehat{f}_2&=-y \cos z \,
\partial_x-y \sin z \, \partial_y+\cos z \, \partial_z.
\end{align}\el
\begin{proof} It is equivalent to show that the rotation defined in the proof of  Proposition \ref{t:rotaz} is $\theta(x,y,z)=z$. The dual basis to our frame $\{f_{1},f_{2},f_{0}\}$ is given by
$$
 \nu_1=-\frac{1}{y}dy, \qquad
 \nu_2=-\frac{1}{y}dx, \qquad
 \nu_0=-\frac{1}{y}dx-dz.
$$
Moreover we have $[f_{1},f_{0}]=[f_{2},f_{0}]=0$ and $[f_{2},f_{1}]=f_{2}+f_{0}$ so that, in equation \eqref{eq:unoforma} we get $c=0,\alpha_1=0,\alpha_2=1$. Hence
  $$d\theta=-\nu_0+\nu_2=dz.$$
\end{proof}

Now we have two canonical frames $\{\widehat{f_1},\widehat{f_2},f_0\}$ and $\{g_1,g_2,g_0\}$, whose Lie algebras satisfy the same commutator relations:
\begin{align}\label{eq:relazioni}
    [\widehat{f}_1,f_0]&=- \widehat{f}_2,    &     &[g_1,g_0]=-g_2, \notag\\
    [\widehat{f}_2,f_0]&=\widehat{f}_1, &   &[g_2,g_0]=g_1,\\
    [\widehat{f}_2,\widehat{f}_1]&= f_0,  & &[g_2,g_1]=0. \notag
\end{align}
Let us consider the two control systems
\begin{align*}
    \dot{q}&=u_1\widehat{f}_1(q)+u_2\widehat{f}_2(q)+u_0f_0(q), \quad
    q \in A^+(\mb{R})\oplus \mb{R}, \\
    \dot{x}&=u_1g_1(x)+u_2g_2(x)+u_0g_0(x), \quad x \in \widetilde{SL}(2).
\end{align*}
and denote with $x_u(t),q_u(t), \ t\in [0,T]$ the solutions of the
equations relative to the same control $u=(u_{1},u_{2},u_{0})$.
Nagano Principle (see \cite{agrcontrol} and also \cite{nagano,sussmannlie,sussmannbook})
ensure that the map \bqn \label{eq:psi}
\til{\Psi}:A^+(\mb{R})\oplus \mb{R} \rightarrow
\widetilde{SL}(2),\qquad q_u(T) \mapsto x_u(T). \eqn that sends the
final point of the first system to the final point of the second
one, is well-defined and does not depend on the control $u$.

Thus we can find the endpoint map of both systems relative to
constant controls, i.e. considering maps \begin{align}
\widetilde{F}&:\R^3 \to A^+(\R)\oplus \R, &   &(t_1,t_2,t_0) \mapsto e^{t_0f_0} \circ e^{t_2\widehat{f}_2} \circ
e^{t_1\widehat{f}_1}(1_A),  \label{eq:traiett1}\\
\widetilde{G}&:\R^3 \to SL(2),     &          &(t_1,t_2,t_0)\mapsto e^{t_0g_0} \circ e^{t_2g_2} \circ e^{t_1g_1}(1_{SL}).
\label{eq:traiett2}
\end{align}
where we denote with
$1_A$ and $1_{SL}$ identity element of $A^+(\R)\oplus \R$ and
$\widetilde{SL}(2)$, respectively.

The composition of these two maps makes the following diagram commutative
\begin{equation}\label{eq:diagriv2}
\xymatrix{
 A^+(\mb{R})\oplus \mb{R} \ar[r]^{\widetilde{\Psi}}  \ar[dr]^{\Psi}  \ar[d]^{\til{F}^{-1}}
&   \widetilde{SL}(2) \ar[d]^{\pi}\\
\mb{R}^{3} \ar[r]^{\til{G}} & SL(2) }
\end{equation}
where $\pi: \widetilde{SL}(2) \to SL(2)$ is the canonical projection
and we set $ \Psi := \pi \circ \widetilde{\Psi}$.

To simplify computation we introduce the rescaled maps
$$F(t):=\til{F}(2t), \qquad G(t):=\til{G}(2t), \qquad t=(t_1,t_2,t_0), $$
and solving differential equations we get from
\eqref{eq:traiett1} the following expressions

\bqn \label{eq:F}
F(t_1,t_2,t_0)=\left(
2e^{-2t_1}\frac{\tanh t_2}{1+\tanh^2 t_2 },\
-e^{-2t_1}\frac{1-\tanh^2t_2}{1+\tanh^2 t_2 },\
2(\arctan(\tanh t_2)-t_0) \right).
\eqn
%To compute $\Psi=G\circ F^{-1}$ we need the inverse function of \eqref{eq:F}, that is%$$
\noindent
The function $F$ is globally invertible on its image and its inverse

$$
F^{-1}(x,y,z)=\left(
-\frac{1}{2} \log \sqrt{x^2+y^2}, \
\tx{arctanh}(\frac{y+\sqrt{x^2+y^2}}{x}), \
\tx{arctan} (\frac{y+\sqrt{x^2+y^2}}{x})-\frac{z}{2}
\right).$$
is defined for every $y<0$ and for every $x$ (it is extended by continuity at $x=0$).

On the other hand, the map  \eqref{eq:traiett2} can be expressed by the product of exponential matrices as follows\ffoot{since we consider left-invariant system, we must multiply matrices on the right.}
\begin{align}\label{eq:G}
G(t_1,t_2,t_0)
=\begin{pmatrix}
e^{t_1} & 0 \\
0 & e^{-t_2} \\
\end{pmatrix}
\begin{pmatrix}
\cosh t_2 & \sinh t_2\\
\sinh t_2 & \cosh t_2\\
\end{pmatrix}
\begin{pmatrix}
\cos t_0 & -\sin t_0 \\
\sin t_0 & \cos t_0 \\
\end{pmatrix}.
%=\begin{pmatrix}
%e^{t_1}(\cos t_0 \cosh t_2 + \sin t_0 \sinh t_2) &
%e^{t_1}(-\sin t_0 \cosh t_2 +\cos t_0 \sinh t_2)\\
%e^{-t_1}(\cos t_0 \sinh t_2 + \sin t_0 \cosh t_2) &
%e^{-t_1}(-\sin t_0 \sinh t_2 +  \cos t_0 \cosh t_2)\\
%\end{pmatrix} %\notag
\end{align}

To simplify the computations, we consider standard polar coordinates $(\rho,\theta)$ on the half-plane $\{(x,y),y<0\}$, where $-\pi/2<\theta<\pi/2$ is the angle that the point $(x,y)$ defines with $y$-axis. In particular, it is easy to see that the expression that appear in $F^{-1}$ is naturally related to these coordinates:
$$\xi=\xi(\theta):=\tan \frac{\theta}{2}=
\begin{cases}
\dfrac{y+\sqrt{x^2+y^2}}{x},~~~~~~ \text{if} ~~~x\neq 0,\\
0,~~~~~~~~~~~~~~~~~~~~~~~\text{if}~~~ x=0.
\end{cases}$$
Hence we can rewrite
$$
F^{-1}(\rho,\theta,z)=\left( -\frac{1}{2} \log \rho, \
\tx{arctanh}\, \xi ,\ \tx{arctan}\, \xi-\frac{z}{2} \right). $$ and
compute the composition $\Psi=G\circ F^{-1} : A^+(\mb{R})\oplus
\mb{R} \longrightarrow SL(2)$. Once we substitute these
expressions in \eqref{eq:G}, the third factor is a rotation matrix
by an angle $\tx{arctan}\, \xi-z/2$. Splitting
this matrix in two consecutive rotations and using standard
trigonometric identities 
$\cos (\tx{arctan} \, \xi)=
\frac{1}{\sqrt{1+\xi^{2}}}, \ \sin (\tx{arctan} \, \xi)=
\frac{\xi}{\sqrt{1+\xi^{2}}}, \ \cosh (\tx{arctanh} \, \xi)=
\frac{1}{\sqrt{1-\xi^{2}}},\ \sinh (\tx{arctanh} \, \xi)=
\frac{\xi}{\sqrt{1-\xi^{2}}}$, for $\xi\in(-1,1)$, 
%
%\begin{gather}
%\cos (\tx{arctan} \, \xi)=
%\frac{1}{\sqrt{1+\xi^{2}}}, \ \sin (\tx{arctan} \, \xi)=
%\frac{\xi}{\sqrt{1+\xi^{2}}}\\
% \ \cosh (\tx{artanh} \, \xi)=
%\frac{1}{\sqrt{1-\xi^{2}}},\ \sinh (\tx{artanh} \, \xi)=
%\frac{\xi}{\sqrt{1-\xi^{2}}}, 
%%\ \xi\in(-1,1)
%\end{gather}
we obtain:
\begin{align*}\label{eq:comp1}
&\Psi(\rho, \theta ,z)=\\
&=\begin{pmatrix}
\rho^{-1/2} & 0 \\
0 & \rho^{1/2} \\
\end{pmatrix}
\begin{pmatrix}
\dfrac{1}{\sqrt{1-\xi^2}} & \dfrac{\xi}{\sqrt{1-\xi^2}}\\
\dfrac{\xi}{\sqrt{1-\xi^2}} & \dfrac{1}{\sqrt{1-\xi^2}}\\
\end{pmatrix}
\begin{pmatrix}
\dfrac{1}{\sqrt{1+\xi^2}} & -\dfrac{\xi}{\sqrt{1+\xi^2}}\\
\dfrac{\xi}{\sqrt{1+\xi^2}} & \dfrac{1}{\sqrt{1+\xi^2}}\\
\end{pmatrix}
\begin{pmatrix}
\cos \dfrac{z}{2} & \sin \dfrac{z}{2} \\[0.3cm]
-\sin \dfrac{z}{2} & \cos \dfrac{z}{2} \\
\end{pmatrix}.
\end{align*}
Then using identities: $\cos \theta=  \dfrac{1-\xi^2}{1+\xi^2}, \ \sin \theta=\dfrac{2\xi}{1+\xi^2}$, we get
\begin{align*}
\Psi(\rho, \theta ,z)&=\begin{pmatrix}
\rho^{-1/2} & 0 \\
0 & \rho^{1/2} \\
\end{pmatrix}
\begin{pmatrix}
\dfrac{1+\xi^2}{\sqrt{1-\xi^4}} & 0\\
\dfrac{2\xi}{\sqrt{1-\xi^4}} & \dfrac{1-\xi^2}{\sqrt{1-\xi^4}}\\
\end{pmatrix}
\begin{pmatrix}
\cos \dfrac{z}{2} & \sin \dfrac{z}{2} \\[0.3cm]
-\sin \dfrac{z}{2} & \cos \dfrac{z}{2} \\
\end{pmatrix}\\
\\
&=\sqrt{\frac{1+\xi^{2}}{1-\xi^{2}}}
\begin{pmatrix}
\rho^{-1/2} & 0 \\
0 & \rho^{1/2} \\
\end{pmatrix}
\begin{pmatrix}
1& 0\\
\dfrac{2\xi}{1+\xi^2} & \dfrac{1-\xi^2}{1+\xi^2}\\
\end{pmatrix}
\begin{pmatrix}
\cos \dfrac{z}{2} & \sin \dfrac{z}{2} \\[0.3cm]
-\sin \dfrac{z}{2} & \cos \dfrac{z}{2} \\
\end{pmatrix}\\
\\
&=\frac{1}{\sqrt{\rho \cos \theta}}
\begin{pmatrix}
1 & 0 \\
0 & \rho \\
\end{pmatrix}
\begin{pmatrix}
1 & 0\\
\sin \theta & \cos \theta\\
\end{pmatrix}
\begin{pmatrix}
\cos \dfrac{z}{2} & \sin \dfrac{z}{2} \\[0.3cm]
-\sin \dfrac{z}{2} & \cos \dfrac{z}{2} \\
\end{pmatrix}\\
\\
&=\frac{1}{\sqrt{\rho \cos\theta}}\begin{pmatrix}
\displaystyle{\cos \frac{z}{2} } &
\displaystyle{\sin \frac{z}{2} }\\[0.2cm]
\displaystyle{\rho \sin (\theta - \dfrac{z}{2})}&
\displaystyle{\rho \cos (\theta - \dfrac{z}{2})}\\
\end{pmatrix}.
\end{align*}

\bl\label{l:q} The set $\Psi^{-1}(I)$ is a normal subgroup of $A^{+}(\R)\oplus \R$.
\el
\begin{proof}
It is easy to show that $\Psi^{-1}(I)=\{F(0,0,2k\pi), k\in \mb{Z}\}$. From \eqref{eq:F} we see that  $F(0,0,2k\pi)=(0,-1,-4k\pi)$ and \eqref{eq:grupL2+R} implies that this is a normal subgroup. Indeed it is enoough to prove that $\Psi^{-1}(I)$ is a subgroup of the centre, that follows from the identity
\begin{equation*}
\begin{pmatrix}
1 & 0 & 0 \\
0 & 1 & 4k\pi \\
0 & 0 & 1 \\
\end{pmatrix}
\begin{pmatrix}
-y & 0 & x\\
0 & 1 & z\\
0 & 0 & 1 \\
\end{pmatrix}=
\begin{pmatrix}
-y & 0 & x\\
0 & 1 & z+4k\pi\\
0 & 0 & 1 \\
\end{pmatrix}=\begin{pmatrix}
-y & 0 & x \\
0 & 1 & z \\
0 & 0 & 1 \\
\end{pmatrix}
\begin{pmatrix}
1 & 0 & 0\\
0 & 1 & 4k\pi\\
0 & 0 & 1 \\
\end{pmatrix}.
\end{equation*}

\end{proof}
\brem With a standard topological argument it is possible to prove that actually $\Psi^{-1}(A)$ is a discrete countable set for every $A\in SL(2)$, and $\Psi$ is a representation of $A^{+}(\R)\oplus \R$ as universal covering of $SL(2)$.
\erem

By Lemma \ref{l:q} the map $\Psi$ is well defined isomorphism between the quotient
$$\frac{A^{+}(\R)\oplus \R}{\Psi^{-1}(I)}\ \simeq \ A^{+}(\R)\times S^{1},$$
and the group $SL(2)$,
defined by restriction of $\Psi$ on $z\in [-2\pi,2\pi]$.

If we consider the new variable $\ph= z/2$, defined on $[-\pi,\pi]$, we can finally write the global isometry as
\begin{gather}\label{eq:uu}
\Psi(\rho,\theta,\ph)=\frac{1}{\sqrt{\rho \cos\theta}}\begin{pmatrix}
\displaystyle{\cos \ph} &
\displaystyle{\sin \ph}\\
\displaystyle{\rho \sin (\theta - \ph)}&
\displaystyle{\rho \cos (\theta - \ph)}\\
\end{pmatrix},
\end{gather} where $(\rho, \theta)\in A^+(\R)$ and $\ph \in S^{1}$.

\brem In the coordinate set defined above we have that $1_A=(1,0,0)$ and
$$\Psi(1_{A})=\Psi(1,0,0)=\begin{pmatrix}
1&
0\\
0&
1\\
\end{pmatrix}=1_{SL}.$$
%By definition $\Psi$ transforms integral curves of the orthonormal frame $\widehat{f}_1$,$\widehat{f}_2$ into integral curves of $g_1,g_2$. Hence $\Psi_{*}(\widehat{f}_{i})=g_{i}$ for $i=1,2$. 
On the other hand $\Psi$ is not a homomorphism since in $A^{+}(\R)\oplus \R$ it holds 
$$\big(\frac{\sqrt{2}}{2},\frac{\pi}{4},\pi\big)\big(\frac{\sqrt{2}}{2},-\frac{\pi}{4},-\pi\big)=1_A,$$ while it can be easily checked from \eqref{eq:uu} that
$$\Psi\big(\frac{\sqrt{2}}{2},\frac{\pi}{4},\pi\big)\Psi\big(\frac{\sqrt{2}}{2},-\frac{\pi}{4},-\pi\big)=
\begin{pmatrix}
2&
0\\
1/2&
1/2\\
\end{pmatrix}
\neq 1_{SL}.$$
\erem

\medskip
\noindent

{\small
\bibliography{bib}

\begin{thebibliography}{10}

\bibitem{miovolume}
A.~Agrachev, D.~Barilari, and U.~Boscain.
\newblock On the {H}ausdorff volume in sub-{R}iemannian geometry.
\newblock {\em Calculus of Variations and Partial Differential Equations},
  pages 1--34, 2011.
\newblock 10.1007/s00526-011-0414-y.

\bibitem{nostrolibro}
A.~Agrachev, D.~Barilari, and U.~Boscain.
\newblock Introduction to {R}iemannian and sub-{R}iemannian geometry ({L}ecture
  {N}otes), http://people.sissa.it/agrachev/agrachev\_files/notes.html.
\newblock 2012.

\bibitem{laplacian}
A.~Agrachev, U.~Boscain, J.-P. Gauthier, and F.~Rossi.
\newblock The intrinsic hypoelliptic {L}aplacian and its heat kernel on
  unimodular {L}ie groups.
\newblock {\em J. Funct. Anal.}, 256(8):2621--2655, 2009.

\bibitem{AAAICM}
A.~A. Agrachev.
\newblock Methods of control theory in nonholonomic geometry.
\newblock In {\em Proceedings of the {I}nternational {C}ongress of
  {M}athematicians, {V}ol.\ 1, 2 ({Z}\"urich, 1994)}, pages 1473--1483.
  Birkh\"auser, Basel, 1995.

\bibitem{agrexp}
A.~A. Agrachev.
\newblock Exponential mappings for contact sub-{R}iemannian structures.
\newblock {\em J. Dynam. Control Systems}, 2(3):321--358, 1996.

\bibitem{localagr}
A.~A. Agrachev, G.~Charlot, J.~P.~A. Gauthier, and V.~M. Zakalyukin.
\newblock On sub-{R}iemannian caustics and wave fronts for contact
  distributions in the three-space.
\newblock {\em J. Dynam. Control Systems}, 6(3):365--395, 2000.

\bibitem{agrcontrol}
A.~A. Agrachev and Y.~L. Sachkov.
\newblock {\em Control theory from the geometric viewpoint}, volume~87 of {\em
  Encyclopaedia of Mathematical Sciences}.
\newblock Springer-Verlag, Berlin, 2004.
\newblock Control Theory and Optimization, II.

\bibitem{bellaiche}
A.~Bella{\"{\i}}che.
\newblock The tangent space in sub-{R}iemannian geometry.
\newblock {\em J. Math. Sci. (New York)}, 83(4):461--476, 1997.
\newblock Dynamical systems, 3.

\bibitem{boscainrossi}
U.~Boscain and F.~Rossi.
\newblock Invariant {C}arnot-{C}aratheodory metrics on {$S^3,\ {\rm SO}(3),\
  {\rm SL}(2)$}, and lens spaces.
\newblock {\em SIAM J. Control Optim.}, 47(4):1851--1878, 2008.

\bibitem{cartan2}
{\'E}.~Cartan.
\newblock Sur la g\'eom\'etrie pseudo-conforme des hypersurfaces de l'espace de
  deux variables complexes {II}.
\newblock {\em Ann. Scuola Norm. Sup. Pisa Cl. Sci. (2)}, 1(4):333--354, 1932.

\bibitem{cartan1}
{\'E}.~Cartan.
\newblock Sur la g\'eom\'etrie pseudo-conforme des hypersurfaces de l'espace de
  deux variables complexes.
\newblock {\em Ann. Mat. Pura Appl.}, 11(1):17--90, 1933.

\bibitem{chow}
W.-L. Chow.
\newblock \"{U}ber {S}ysteme von linearen partiellen {D}ifferentialgleichungen
  erster {O}rdnung.
\newblock {\em Math. Ann.}, 117:98--105, 1939.

\bibitem{falbel}
E.~Falbel and C.~Gorodski.
\newblock Sub-{R}iemannian homogeneous spaces in dimensions {$3$} and {$4$}.
\newblock {\em Geom. Dedicata}, 62(3):227--252, 1996.

\bibitem{gersh}
V.~Gershkovich and A.~Vershik.
\newblock Nonholonomic manifolds and nilpotent analysis.
\newblock {\em J. Geom. Phys.}, 5(3):407--452, 1988.

\bibitem{gromov}
M.~Gromov.
\newblock Carnot-{C}arath\'eodory spaces seen from within.
\newblock In {\em Sub-{R}iemannian geometry}, volume 144 of {\em Progr. Math.},
  pages 79--323. Birkh\"auser, Basel, 1996.

\bibitem{jacob}
N.~Jacobson.
\newblock {\em Lie algebras}.
\newblock Interscience Tracts in Pure and Applied Mathematics, No. 10.
  Interscience Publishers (a division of John Wiley \& Sons), New York-London,
  1962.

\bibitem{sachkovmois}
I.~Moiseev and Y.~L. Sachkov.
\newblock Maxwell strata in sub-{R}iemannian problem on the group of motions of
  a plane.
\newblock {\em ESAIM Control Optim. Calc. Var.}, 16:380--399, 2010.

\bibitem{montgomery}
R.~Montgomery.
\newblock {\em A tour of subriemannian geometries, their geodesics and
  applications}, volume~91 of {\em Mathematical Surveys and Monographs}.
\newblock American Mathematical Society, Providence, RI, 2002.

\bibitem{nagano}
T.~Nagano.
\newblock Linear differential systems with singularities and an application to
  transitive {L}ie algebras.
\newblock {\em J. Math. Soc. Japan}, 18:398--404, 1966.

\bibitem{rashevsky}
P.~Rashevsky.
\newblock Any two points of a totally nonholonomic space may be connected by an
  admissible line.
\newblock {\em Uch. Zap. Ped Inst. im. Liebknechta}, 2:83--84, 1938.

\bibitem{sachkov}
Y.~L. Sachkov.
\newblock Conjugate and cut time in the sub-{R}iemannian problem on the group
  of motions of a plane.
\newblock {\em ESAIM Control Optim. Calc. Var.}, 16:1018--1039, 2010.

\bibitem{stric}
R.~S. Strichartz.
\newblock Sub-{R}iemannian geometry.
\newblock {\em J. Differential Geom.}, 24(2):221--263, 1986.

\bibitem{stric2}
R.~S. Strichartz.
\newblock Corrections to: ``{S}ub-{R}iemannian geometry''.
\newblock {\em J. Differential Geom.}, 30(2):595--596, 1989.

\bibitem{sussmannlie}
H.~J. Sussmann.
\newblock An extension of a theorem of {N}agano on transitive {L}ie algebras.
\newblock {\em Proc. Amer. Math. Soc.}, 45:349--356, 1974.

\bibitem{sussmannbook}
H.~J. Sussmann.
\newblock Lie brackets, real analyticity and geometric control.
\newblock In {\em Differential geometric control theory ({H}oughton, {M}ich.,
  1982)}, volume~27 of {\em Progr. Math.}, pages 1--116. Birkh\"auser Boston,
  Boston, MA, 1983.

\end{thebibliography}
\bibliographystyle{abbrv}
}

\end{document}